\documentclass[12pt]{amsart}
\usepackage{amsmath, amssymb, latexsym, tikz, tikz-cd, amsthm, mathrsfs, color, mathtools, soul, float}
\usepackage[hyphens]{url}
\usepackage[T1]{fontenc}
\usepackage{lmodern}

\usepackage[top=2in, bottom=1.5in, left=1in, right=1in]{geometry}

\newcommand{\nc}{\newcommand}

\nc{\thusfar}{\par\bigskip\centerline{\my{--- Edited thus far ---}}\par\bigskip}
\nc{\lei}{\le^\oo}
\nc{\card}[1]{\left|#1\right|}
\nc{\medcard}[1]{\biggl|\,#1\,\biggr|}
\nc{\smallcard}[1]{|\,#1\,|}
\nc{\bbN}{\mathbb{N}}
\nc{\beq}{\begin{eqnarray*}}\nc{\eeq}{\end{eqnarray*}}
\nc{\mbq}{\mb{?}}
\nc{\mb}[1]{{\mbox{\textbf{#1}}}}
\nc{\nop}{$\times$}
\nc{\fbn}{\!\!\fbox{\!\nop\!}\!\!}
\nc{\yup}{\checkmark}
\nc{\forces}{\Vdash}
\nc{\name}[1]{\dot{#1}}
\nc{\tf}{\my{FINISHED THUS FAR}}
\nc{\FU}{Fr\'echet--Urysohn}
\nc{\gs}{$\gamma$~space}
\nc{\Ga}{\Gamma}\nc{\Om}{\Omega}
\nc{\smallbinom}[2]{\begin{psmallmatrix} #1\\ #2 \end{psmallmatrix}}
\nc{\bgamma}{\smallbinom{\Om}{\Ga}}

\nc{\productive}[2]{\bigl(#1,\allowbreak #2\bigr)^\x}
\nc{\Sel}{\mathsf{S}}
\nc{\sset}[2]{\{\,#1 : #2\,\}}
\nc{\smb}[1]{{\!\!\mb{#1}\!\!}}
\nc{\medset}[2]{{\biggl\{\,#1 : #2\,\biggr\}}}
\nc{\smallmedset}[2]{{\bigl\{\,#1 : #2\,\bigr\}}}
\nc{\set}[2]{{\left\{\,#1 : #2\,\right\}}}
\nc{\seq}[2]{{\la\, #1 : #2\,\ra}}
\nc{\eseq}[1]{#1_0, \allowbreak #1_1, \allowbreak\dots} 
\nc{\cube}{(\Cantor)^\bbN}
\nc{\Match}{\op{Match}}
\nc{\concat}[1]{\hat{\phantom{a}}\langle #1\rangle}
\nc{\poset}{\mathbb{P}}
\nc{\fn}[1]{{\op{Fn}(#1\times\w,2)}}
\nc{\linadd}{\op{linadd}}
\nc{\nonprod}{\non^\x}
\nc{\alephes}{{\aleph_0}}
\nc{\my}[1]{\marginpar{\textcolor{red}{***}}\textcolor{red}{#1}}
\nc{\later}[1]{{\color{green} #1}}
\nc{\BTs}[1]{{\color{green} #1 (BT)}}
\nc{\Cp}{\op{C}_\mathrm{p}}
\nc{\Bp}{\op{B}_p}
\nc{\Pa}[8]{\bibitem{#1} {#2}, \emph{#3}, {#4} \textbf{#5} ({#6}), {#7}--{#8}.}
\nc{\tPa}[5]{\bibitem{#1} {#2}, \emph{#3}, {#4}, to appear.}
\nc{\sPa}[4]{\bibitem{#1} {#2}, \emph{#3}, {#4}, submitted.}
\nc{\Bc}[9]{\bibitem{#1} {#2}, \emph{#3}, in: \textbf{#4} (#5), #6 #7, #8--#9.}
\nc{\fD}{\mathfrak{D}}
\nc{\fX}{\mathfrak{X}}
\nc{\Onbd}{\Op_{\mathrm{nbd}}} 
\nc{\Omnb}{\Om_{\mathrm{nbd}}} 
\nc{\od}{\mathfrak{od}}
\nc{\Setting}[7]{\xymatrix@R=4pt@C=7pt{#1\ar@{-}[r]&#2\ar@{-}[r]&#3\\&#4\ar@{-}[u]\\
		#5\ar@{-}[uu]\ar@{-}[r] & #6\ar@{-}[u]\ar@{-}[r] & #7\ar@{-}[uu]}}
\nc{\mx}[1]{\begin{matrix}#1\end{matrix}}
\nc{\plim}{p\txt{-}\lim}
\nc{\Bgp}{{\Z^\bbN}}
\nc{\Cgp}{{{\Z_2}^\bbN}}
\nc{\Cite}[1]{\textbf{[#1]}}
\nc{\Next}[1]{{#1^+}}
\nc{\cFin}{\mathrm{cF}}
\nc{\intvl}[2]{{[#1(#2),\allowbreak #1(#2\!+\!1))}}
\nc{\Bdd}{\mathbf{B}}
\nc{\Dfin}{\mathfrak{D}_\mathrm{fin}}
\nc{\grbl}{{\mbox{\textit{\tiny gp}}}}
\nc{\bbP}{\mathbb{P}}
\nc{\BOfat}{\B_{\Om_{\mathrm{fat}}}}
\nc{\Bgood}{\B_{\mathrm{good}}}
\nc{\compactN}{\cl{\mathbb{N}}}
\nc{\blocks}[2]{\op{cl}_{#2}(#1)}
\nc{\blocksplus}[2]{\op{cl}^+_{#2}(#1)}
\nc{\arx}[1]{\texttt{http://arxiv.org/math/#1}}
\nc{\bq}{\begin{quote}}
\nc{\eq}{\end{quote}}
\nc{\cl}[1]{\overline{#1}}
\nc{\CH}{the Continuum Hypothesis}
\nc{\MA}{Martin's Axiom}
\nc{\Bfat}{\B_\mathrm{fat}}
\nc{\inv}{^{-1}}
\nc{\Cantor}{{2^\w}}
\nc{\bP}{\mathbf{P}}
\nc{\bof}{\op{\fb}}
\nc{\dof}{\op{\fd}}
\nc{\bofF}{\bof(\cF)}
\nc{\sr}[3]{\underset{\mbox{#3}}{\mbox{#1}}}
\nc{\gp}{\binom{\Om}{\Ga}}
\nc{\gpsmall}{\mbox{$\gp$}}
\nc{\gig}{\gimel}
\nc{\gns}{\sone(\Om,\gig)}
\nc{\nsr}[2]{#1}
\nc{\Srg}{{\mathbb{S}}}
\nc{\Srgs}{{\mathbb{S}^*}}
\nc{\NN}{{\w^{\w}}}
\nc{\ZN}{{\Z^{\bbN}}}
\nc{\NNup}{{\bbN^{\uparrow\bbN}}}
\nc{\Pof}{\op{P}}
\nc{\PN}{{\Pof(\w)}}
\nc{\roth}{{[\w]^{\w}}}
\nc{\Fin}{\mathrm{Fin}} 
\nc{\ici}{[\bbN]^{ \infty, \infty}}
\nc{\Inc}{{\compactN^{\uparrow\bbN}}}
\nc{\powInc}[1]{{\big(\Inc\big)^{#1}}}
\nc{\powFin}[1]{{\big(\Fin\big)^{#1}}}
\nc{\powPN}[1]{{\big(\PN\big)^{#1}}}
\nc{\NcompactN}{{\compactN^\bbN}}
\nc{\Uarrow}{\smash{\big\uparrow}}
\nc{\LE}{\preccurlyeq}
\nc{\GE}{\succcurlyeq}
\nc{\op}{\operatorname}
\nc{\im}{\op{im}}
\nc{\Span}{\op{span}}
\nc{\maxfin}{\op{maxfin}}
\nc{\ran}{\op{range}}
\nc{\iso}{\cong}
\nc{\Madd}{{\M}^\star}
\nc{\cI}{\mathcal{I}}
\nc{\cJ}{\mathcal{J}}
\nc{\scrA}{\mathscr{A}}
\nc{\scrB}{\mathscr{B}}
\nc{\scrC}{\mathscr{C}}
\nc{\scrD}{\mathscr{D}}
\nc{\scrF}{\mathscr{F}}
\nc{\scrK}{\mathscr{K}}
\nc{\A}{\forall}
\nc{\B}{\mathrm{B}}
\nc{\cB}{\mathcal{B}}
\nc{\bB}{\mathbf{B}}
\nc{\BS}{\mathbf{B}(\mathcal{S})}
\nc{\BF}{\mathbf{B}(\mathcal{F})}
\nc{\BU}{\mathbf{B}(\mathcal{U})}
\nc{\cSp}{\mathcal{S}^+}
\nc{\cFp}{\mathcal{F}^+}
\nc{\cUp}{\mathcal{U}^+}

\nc{\BG}{\B_\Ga}
\nc{\BL}{\B_\Lambda}
\nc{\BT}{\B_\Tau}
\nc{\BTstar}{\B_{\Tau^*}}
\nc{\BO}{\B_\Om}
\nc{\DO}{\cD_\Om}
\nc{\KO}{\cK_\Om}
\nc{\CG}{C_\Ga}
\nc{\CL}{C_\Lambda}
\nc{\CT}{C_\Tau}
\nc{\CTstar}{C_{\Tau^*}}
\nc{\CO}{C_\Om}
\nc{\COgp}{C_{\Om^{\grbl}}}
\nc{\CLgp}{C_{\Lambda^{\grbl}}}
\nc{\BOgp}{\B_{\Om}^{\grbl}}
\nc{\BLgp}{\B_{\Lambda^{\grbl}}}
\nc{\sfC}{\mathsf{C}}
\nc{\sfD}{\mathsf{D}}
\nc{\bD}{\mathbf{D}}
\nc{\Tau}{\mathrm{T}}
\nc{\cA}{\mathcal{A}}
\nc{\cK}{\mathcal{K}}
\nc{\cD}{\mathcal{D}}
\nc{\cF}{\mathcal{F}}
\nc{\cS}{\mathcal{S}}
\nc{\cT}{\mathcal{T}}
\nc{\cG}{\mathcal{G}}
\nc{\cY}{\mathcal{Y}}
\nc{\J}{\mathcal{J}}
\nc{\cL}{\mathcal{L}}
\nc{\cM}{\mathcal{M}}
\nc{\cN}{\mathcal{N}}
\nc{\cH}{\mathcal{H}}
\nc{\cO}{\mathcal{O}}
\nc{\Op}{\mathrm{O}}
\nc{\rmA}{\mathrm{A}}
\nc{\rmF}{\mathrm{F}}
\nc{\rmB}{\mathrm{B}}
\nc{\rmD}{\mathrm{D}}
\nc{\rmP}{\mathrm{P}}
\nc{\cC}{\mathcal{C}}
\nc{\cP}{\mathcal{P}}
\nc{\bbQ}{\mathbb{Q}}
\nc{\bbR}{\mathbb{R}}
\nc{\cU}{\mathcal{U}}
\nc{\cQ}{\mathcal{Q}}
\nc{\bbC}{\mathbb{C}}
\nc{\Un}{\bigcup}
\nc{\cV}{\mathcal{V}}
\nc{\cW}{\mathcal{W}}
\nc{\Z}{{\mathbb Z}}
\nc{\Impl}{\Rightarrow}
\long\def\forget#1\forgotten{\marginpar{\textcolor{green}{Forgetting...}}}
\nc{\ft}{\mathfrak{t}}
\nc{\fb}{\mathfrak{b}}
\nc{\fc}{\mathfrak{c}}
\nc{\fd}{\mathfrak{d}}
\nc{\fg}{\mathfrak{g}}
\nc{\oo}{\infty}
\nc{\fr}{\mathfrak{r}}
\nc{\fk}{\mathfrak{k}}
\nc{\bidi}{\mathfrak{bidi}}
\nc{\fu}{\mathfrak{u}}
\nc{\fh}{\mathfrak{h}}
\nc{\fp}{\mathfrak{p}}
\nc{\fj}{\mathfrak{j}}
\nc{\fs}{\mathfrak{s}}
\nc{\w}{\omega}
\nc{\x}{\times}
\nc{\Iff}{\Leftrightarrow}

\nc{\nin}{\notin}
\nc{\cat}{\hat{\ }}
\nc{\sub}{\subseteq}
\nc{\spst}{\supseteq}
\nc{\sm}{\setminus}
\nc{\as}{\subseteq^*}
\nc{\les}{\le^*}
\nc{\leinf}{\le^{\infty}}
\nc{\leS}{\le_S}
\nc{\leF}{\le_{\mathcal{F}}}
\nc{\leU}{\le_{U}}
\nc{\rest}{\restriction}

\nc{\la}{\langle}
\nc{\ra}{\rangle}
\nc{\E}{\exists}
\nc{\cov}{\op{cov}}
\nc{\add}{\op{add}}
\nc{\cof}{\op{cof}}
\nc{\cf}{\op{cf}}
\nc{\non}{\op{non}}
\nc{\unif}{\op{non}}
\nc{\COV}{\op{COV}}
\nc{\ADD}{\op{ADD}}
\nc{\COF}{\op{COF}}
\nc{\NON}{\op{NON}}
\nc{\impl}{\to}
\nc{\Lp}{\mathcal{L_\p}}
\nc{\Wlog}{without loss of generality}
\newtheorem{thm}{Theorem}[section]
\nc{\bthm}{\begin{thm}} \nc{\ethm}{\end{thm}}
\newtheorem{prop}[thm]{Proposition}
\nc{\bprp}{\begin{prop}} \nc{\eprp}{\end{prop}}
\newtheorem{fact}[thm]{Fact}
\nc{\bfct}{\begin{fact}} \nc{\efct}{\end{fact}}
\newtheorem{prob}[thm]{Problem}
\nc{\bprb}{\begin{prob}} \nc{\eprb}{\end{prob}}
\newtheorem{lem}[thm]{Lemma}
\nc{\blem}{\begin{lem}} \nc{\elem}{\end{lem}}
\newtheorem{claim}[thm]{Claim}
\nc{\bclm}{\begin{claim}} \nc{\eclm}{\end{claim}}
\newtheorem{cor}[thm]{Corollary}
\nc{\bcor}{\begin{cor}} \nc{\ecor}{\end{cor}}
\newtheorem{conj}[thm]{Conjecture}
\nc{\bcnj}{\begin{conj}} \nc{\ecnj}{\end{conj}}
\theoremstyle{definition}
\newtheorem{defn}[thm]{Definition}
\nc{\bdfn}{\begin{defn}} \nc{\edfn}{\end{defn}}
\newtheorem{obs}[thm]{Observation}
\nc{\bobs}{\begin{obs}} \nc{\eobs}{\end{obs}}
\newtheorem{rem}[thm]{Remark}
\nc{\brem}{\begin{rem}} \nc{\erem}{\end{rem}}
\newtheorem{cnv}[thm]{Convention}
\nc{\bcnv}{\begin{cnv}} \nc{\ecnv}{\end{cnv}}
\newtheorem{exam}[thm]{Example}
\nc{\bexm}{\begin{exam}} \nc{\eexm}{\end{exam}}
\nc{\bpf}{\begin{proof}} \nc{\epf}{\end{proof}}
\nc{\be}{\begin{enumerate}}
\nc{\ee}{\end{enumerate}}
\nc{\bi}{\begin{itemize}}
\nc{\bimy}{\my{\begin{itemize}}
		\nc{\eimy}{\end{itemize}}}
\nc{\itm}{\item}
\nc{\ei}{\end{itemize}}
\nc{\Subsection}[1]{\goodbreak\subsection*{#1}}

\nc{\sone}{\mathsf{S}_1}
\nc{\sfin}{\mathsf{S}_\mathrm{fin}}
\nc{\ufin}{\mathsf{U}_\mathrm{fin}}
\nc{\Split}{\mathrm{Split}}

\nc{\gone}{\mathsf{G}_1}   
\nc{\Succ}{\mathrm{S}} 
\nc{\gfin}{\mathsf{G}_\mathrm{fin}}

\DeclareMathOperator{\eexists}{\exists}
\DeclareMathOperator{\fforall}{\forall}
\nc{\Exists}[1]{\bigl(\eexists #1\bigr)}
\nc{\Forall}[1]{\fforall #1\ }
\nc{\Foralm}[1]{\fforall^* #1\ }
\nc{\End}[1]{#1}
\nc{\supp}{\op{supp}}
\nc{\bfP}{\mathbf{P}}

\nc{\Alice}{{\textsc{Alice}}}
\nc{\Bob}{{\textsc{Bob}}}
\DeclareMathOperator{\M}{M}

\title{Productively Scheepers spaces and their relatives}

\thanks{The research of the third author
	was funded in whole by the Austrian Science Fund (FWF) [10.55776/I5930 and 10.55776/PAT5730424].
	The research of the second author was funded by the National Science Center, Poland Weave-UNISONO call in the Weave programme
	Project: Set-theoretic aspects of topological selections 2021/03/Y/ST1/00122.
}

\author[M. K{\l}ad\'{z}--Duda]{Marta K{\l}ad\'{z}--Duda}
\address{Marta K{\l}ad\'{z}--Duda, Institute of Mathematics, University of Warsaw, Banacha 2, 02-097, Warsaw, Poland}
\email{martakladz3@gmail.com}

\author[P. Szewczak]{Piotr Szewczak}
\address{Piotr Szewczak, 
 Institute of Mathematics, University of Warsaw, Banacha 2, 02-097, Warsaw, Poland
}
\email{p.szewczak@wp.pl}
\urladdr{https://piotrszewczak.pl}

\author[L. Zdomskyy]{Lyubomyr Zdomskyy}
\address{Lyubomyr Zdomskyy, Institut f\"ur Diskrete Mathematik und Geometrie, Technische Universit\"at Wien, Wiedner Hauptstrasse 8-10/104, 1040 Wien, Austria.}
\email{lzdomsky@gmail.com}
\urladdr{https://dmg.tuwien.ac.at/zdomskyy/}

\subjclass[2010]{Primary: 54D20; 
	Secondary: 03E17. 
}

\keywords{Menger, Scheepers, Hurewicz, products}

\begin{document}
	
	\maketitle
	
	\begin{abstract}
		We prove that assuming $\fb=\fd$, in the class of hereditarily Lindel\"of spaces, each productively Scheepers space is productively Hurewicz. 
		This statement remains true for all general topological spaces, assuming that $\mathfrak{d}=\aleph_1$. 
        We also show that if near coherence of filters holds, then the Scheepers property is equivalent to the Menger property parametrized by any ultrafilter.
	\end{abstract}
	
	\section{Introduction}
	
	\subsection{Landscape of combinatorial covering properties}	
	By \emph{space} we mean an infinite Tychonoff topological space.
	Let $\cU$ be a cover of a space $X$.
	The family $\cU$ is a \emph{$\gamma$-cover} of $X$, if it is infinite and each point from $X$ belongs to all but finitely many sets from $\cU$.
	The family $\cU$ is an \emph{$\w$-cover} of $X$, if $X\notin \cU$ and each finite subset of $X$ is contained in a set from $\cU$.
	We consider the following classical combinatorial covering properties.
	A space $X$ is \emph{Hurewicz} (\emph{Scheepers},\emph{ Menger}), if for every sequence $\eseq{\cU}$ of open covers of $X$, none with a finite subcover, there are finite sets $\cF_0\sub \cU_0, \cF_1\sub\cU_1,\dotsc$ such that the family $\sset{\Un\cF_n}{n\in \w}$ is a $\gamma$-cover ($\w$-cover, cover) of $X$.  
	We have the following implications:
	\[
	\sigma\text{-compact } \Longrightarrow \text{ Hurewicz } \Longrightarrow \text{ Scheepers }\Longrightarrow \text{ Menger }
	\]
	and there are examples of sets of reals showing that the first two implications are not reversible in ZFC (Just--Miller--Scheepers--Szeptycki~\cite[Theorem 5.1]{just1995combinatoricsopencoversii}, Bartoszy\'nski--Shelah\cite{BARTOSZYNSKI2001243}, Tsaban--Zdomskyy~\cite[Theorem 3.9]{Tsaban_2008}).
Definitions of the cardinal characteristics of the continuum used in the paper and relations between them, can be found in a work of Blass~\cite{Blass2010}.
Assuming that $\fu<\fg$ (which holds, e.g., in the Miller model) Scheepers and Menger coincide~\cite[Corollary 2]{zdomsky2004semifilterapproachselectionprinciples}; 
on the other hand, assuming that $\fd\leq\fr$, there is a Menger set of reals that is not Scheepers~\cite[Theorem 2.1]{szewczak2020finitepowersproductsmenger}.
	
\subsection{Products}
	A topological property $\bP$ is \emph{productive} if every product of two spaces with the property $\bP$, has the property $\bP$. 
 Throughout the paper we work in the realm of hereditarily Lindel\"of spaces.
	In the Miller model (where Menger is equivalent to Scheepers) Menger is productive but Hurewicz is not ~\cite[Theorem ~1.1]{Zdomskyy_2018},~\cite[Theorem ~1.7]{zdomskyy2019selectionprincipleslavermiller}.
	In the Laver model Hurewicz is productive but Scheepers and Menger are not~\cite[Theorem 1.2]{Repov__2025}.
	If the Continuum Hypothesis holds, then none of the properties Hurewicz, Scheepers, Menger is productive~\cite[Theorem~2.12]{just1995combinatoricsopencoversii}.
	
	Let $\bP$ be a topological property.
	A space $X$ is \emph{productively $\bP$}, if for every space $Y$ with the property $\bP$, the product space $X\x Y$ has the property $\bP$.
	Szewczak--Tsaban proved that assuming $\fb=\fd$ any productively Menger space is productively Hurewicz~\cite[Theorem 4.8]{szewczak2017productsgeneralmengerspaces}. 
	Recently, Repov\v{s}--Zdomskyy proved that each Hurewicz set is productively Menger in the Laver model {~\cite[Theorem 1.1]{Repov__2025}}.
	A combination of these results, delivers the above mentioned result that in the Laver model (where $\fb=\fd$), Hurewicz is productive.
	Repov\v{s}--Zdomskyy also proved that it is consistent with the Continuum Hypothesis that there is a productively Hurewicz space that is not productively Scheepers and not productively Menger~\cite[Theorem 1.2]{Repov__2025}.
	
	\subsection{Main results}
	The aim of this paper is to complete this picture, considering productively Scheepers spaces.
	
	We use the following auxiliary properties which lie (formally) between Hurewicz and Scheepers. 
    By \emph{ultrafilter} we mean a free ultrafilter.
	Let $U$ be an ultrafilter on $\w$.
	A space $X$ is \emph{$U$-Menger}, if for every sequence $\eseq{\cU}$ of open covers of $X$, none with a finite subcover, there are finite sets $\cF_0\sub\cU_0,\cF_1\sub\cU_1,\dotsc$ such that the sets $\sset{n}{x\in\Un\cF_n}$ are in the ultrafilter $U$ for all $x\in X$. 
	\[
	\text{Hurewicz}\Longrightarrow U\text{-Menger}\Longrightarrow\text{Scheepers}\Longrightarrow\text{Menger}.
	\]
	These properties turn out to be extremely useful in constructing counterexamples for productivity of properties, e.g.,  Szewczak--Tsaban--Zdomskyy showed that assuming $\fd\leq\fr$ and regularity of $\fd$, there are two sets which are Menger in all finite powers (in particular they are Scheepers), whose product space is not Menger~\cite[Theorem 2.5]{szewczak2020finitepowersproductsmenger}.
	We will establish some relations between productively Hurewicz, productively $U$-Menger and productively Menger spaces.
	It is consistent\footnote{{The $U$-Menger property is not considered directly in \cite{zdomsky2004semifilterapproachselectionprinciples}, so we refer the reader to Lemma~\ref{expl_nachher} for more details regarding the equivalence of Scheepers and $U$-Menger for some (all) ultrafilters $U$ in models of $\mathfrak u<\mathfrak g$.}} that properties $U$-Menger for all ultrafilters $U$, Scheepers and Menger coincide~\cite[Corollary 2]{zdomsky2004semifilterapproachselectionprinciples} (e.g., under $\fu<\fg$).
	On the other hand, assuming that $\fb=\fd$, there are ultrafilters $U$ and $V$ such that the properties $U$-Menger and $V$-Menger are different (\cite[Theorem~5.3(2)]{Szewczak2017}, \cite[Theorem~3.3]{rezdo}).
	
	\begin{figure}[H]
		\begin{tikzcd}[ampersand replacement=\&,column sep=2.5cm]
			{}\& \begin{matrix}\text{productively}\\ \text{Scheepers}\arrow[ld, "\textbf{Theorem \ref{result}}", swap]\end{matrix}\arrow[d, dashed, "\textbf{\large ?}"]\& {}\\
			\begin{matrix}\text{productively}\\ \text{Hurewicz}\end{matrix}\& \arrow[l, left,"\textbf{Theorem \ref{result}}", swap]\begin{matrix}\text{productively}\\ U\text{-Menger} \\\text{for all }U\end{matrix}\&\arrow[l, "\textbf{Theorem \ref{thm:pMpUM}}", swap]\begin{matrix}\text{productively}\\ \text{Menger} \end{matrix}\arrow[lu, dashed, "\textbf{\large ?}", swap]\arrow[ll, "\textbf{~\cite[Theorem 4.8]{szewczak2017productsgeneralmengerspaces}}", bend left=25]
		\end{tikzcd}
		\caption{Relations in the class of hereditarily Lindel\"of spaces, assuming that $\fb=\fd$. }
	\end{figure}
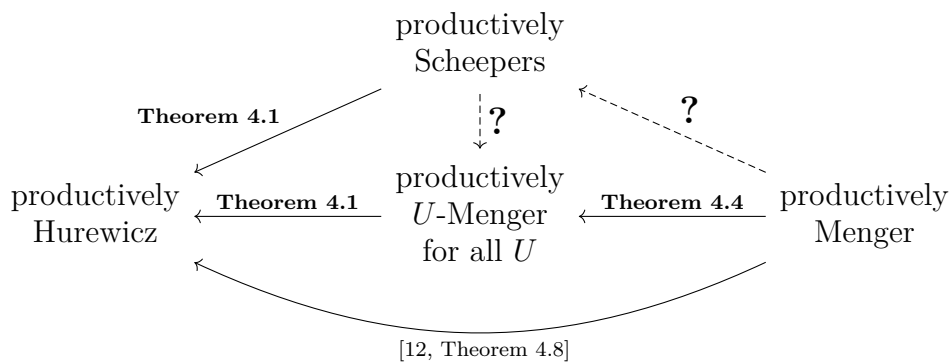

    During our investigations we realized that the Scheepers property is equivalent to the $U$-Menger property for all ultrafilters $U$ under a weaker assumption than $\fu<\fg$, i.e., if \emph{near coherence of filters} holds. We discuss this in Section~\ref{sec:Sch-UMen}.
	
	\section{Combinatorial characterizations of properties}
	
	The properties considered here were succesfully studied using the combinatorial structure of the Baire space.
	Let $\PN$ be the power set of the set of natural numbers $\w$.
	Identifying each set from $\PN$ with its characteristic function, which is an element of the Cantor cube $\Cantor$, we introduce a topology on $\PN$.
	Let $\roth$ be the family of all infinite subsets of $\w$.
	We identify each element of $\roth$ with the increasing enumeration of its elements, an element of the Baire space $\NN$.
	Then $\roth$ can be viewed as a subset of $\PN$ and $\NN$ and the topologies on $\roth$, inherited from $\PN$ and $\NN$, are the same.
	
	For $a,b\in\roth$, we write $a\les b$, if the set $\sset{n}{a(n)\leq b(n)}$ is cofinite.
	A set $A\sub \roth$ is \emph{bounded}, if there is a function $b\in\roth$ such that $a\les b$ for all functions $a\in A$.
	A subset of $\roth$ is \emph{unbounded} if it is not bounded.
	Let $\fb$ be the minimal cardinality of an unbounded set.
	A set $D\sub \roth$ is \emph{dominating}, if for each function $a\in\roth$ there is a function $d\in D$ with $a\les d$.
	Let $\fd$ be the minimal cardinality of a dominating set.
	Let $U$ be an ultrafilter on $\w$. 
	For $a,b\in\roth$ we write $a\leU b$, if the set $\sset{n}{a(n)\leq b(n)}$ is in the ultrafilter $U$.
	A set $A\sub \roth$ is \emph{$\leU$-bounded} in there is a function $b\in \roth$ such that $a\leU b$ for all $a\in A$.
	A subset of $\roth$ is \emph{$\leU$-unbounded} if it is not $\leU$-bounded.
	Let $\bof(U)$ be the minimal cardinality of a $\leU$-unbounded set.
	We have the following inequalities for all ultrafilters $U$:
	\[
	\fb\leq\bof(U)\leq \fd.
	\]
	The numbers $\fb$ and $\fb(U)$ are regular and regularity of $\fd$ is independent from ZFC.
	
	Let $X$ and $Y$ be spaces.
	A set-valued map $\Phi$ that assigns to each element $x\in X$ a compact set $\Phi(x)\sub Y$ is \emph{compact-valued upper semicontinuous} (\emph{cusco} for short) if for every open set $U\sub Y$, the set $\Phi\inv[U]:=\sset{x\in X}{\Phi(x)\sub U}$ is open in $X$.
	
\bthm[Miller, Tsaban, Zdomskyy {\cite[Theorem 7.3, Theorem 7.5]{Miller2014}}]\label{thm:Hur}
Let $X$ be a space.
The following assertions are equivalent.
\begin{enumerate}
\item The space $X$ is Hurewicz ($U$-Menger, Scheepers, Menger).
\item Every cusco image of $X$ in $\roth$ is bounded ($\leU$-bounded, not finitely-dominating, not dominating).
\end{enumerate}
\ethm

	It follows directly from these characterizations that each of the considered properties has a \emph{critical cardinality}, i.e., the minimal cardinality of a subset of $\roth$ that does not have the given property. 
	For Hurewicz it is $\fb$, for $U$-Menger it is $\bof(U)$, and for Scheepers and Menger it is $\fd$.
	Moreover, in ZFC, there is a \emph{nontrivial} Hurewicz, $U$-Menger and  Menger set of reals, i.e., a set of cardinality greater or equal to the critical cardinality of the corresponding property~\cite[Corollary 2.16]{Tsaban_2008}.
	If $\fd$ is regular, then the same holds for the Scheepers property~\cite[Theorem~4.6]{Tsaban_2008}.
    It remains open whether this also holds when $\fd$ is singular.
   
	\section{Combinatorial tools}
	
	A subset of $\roth$ is \emph{$\fd$-unbounded}~\cite[Definition~2.1]{Szewczak2017}, if it has size at least $\fd$ and every bounded subset of this set has size smaller than $\fd$. 
	Let $\Fin$ be the family of all finite subsets of $\w$.
	A $\fd$-unbounded set exists in ZFC~\cite[Lemma~17]{BartoszynskiTsaban2005} and if a set $Y\sub\roth$ is $\fd$-unbounded, then the set $Y\cup\Fin$ is Menger~\cite[Corollary 2.4]{Szewczak2017}.
	
	\bthm[Szewczak--Tsaban~{\cite[Theorem 2.7]{Szewczak2017}}]
	Assume that $\cf(\fd)=\fd$.
	Let $X\sub\roth$ be a set containing a $\fd$-unbounded set.
	Then there exists a $\fd$-unbounded set $Y\sub\roth$, such that the set $X\times (Y\cup\Fin)$ is not Menger.
	\label{Sz-T}
	\ethm
	
	For our purposes we need that the $\fd$-unbounded set $Y$ in Theorem~\ref{Sz-T} has a stronger combinatorial structure.
	Let $U$ be an ultrafilter.
	A set $Y\sub\roth$ is a \emph{$U$-scale}, if for each function $b\in\roth$ the inequality $b\leU y$ holds for all but less than $\fb(U)$ functions $y\in Y$. Such a set exists in ZFC~\cite[Lemma 2.9]{Tsaban_2008}.
	If $Y\sub\roth$ is a $U$-scale then the set $Y\cup\Fin$ is $U$-Menger~\cite[Lemma 2.14]{Tsaban_2008} and even productively $U$-Menger ~\cite[Lemma 4.1]{Tsaban_2008}.
    
    For $a,b\in\roth$ let $[a\leq b]:=\sset{n}{a(n)\leq b(n)}$.
	Let $F\sub\roth$ be a set.
	For $a,b \in\roth$ we write $a\leq_F b$, if $[a\leq b]\in F$.
	For a set $A\sub\roth$ and $b\in \roth$, we write $A\le_F b$, if $a\le_F b$ for all $a\in A$.
	We follow this convention for any binary relation on $\roth$.
	For $a,b\in\roth$ we write $a\leinf b$, if the set $[a\leq b]$ is infinite and $a\leq b$, if $[a\leq b]=\w$.
	
	\bthm\label{thm:U}
	Assume that $\cf\fd=\fd$. 
	Each $\fd$-unbounded set is $\le_U$-unbounded for some ultrafilter $U$ with $\bof(U)=\fd$.
	\label{ultra}
	\ethm

	\bpf
	Let $A\sub \roth$ be a $\fd$-unbounded set
	and $D=\sset{d_\alpha}{\alpha<\fd}$ be a dominating set.
	We may assume that $\card{A}=\fd$.
	
	Recursively, for each $\alpha<\fd$ define elements $a_\alpha\in A$ and sets $F_\alpha\sub\roth$  with the properties:
	\begin{enumerate}
		\item $F_\alpha$ is closed under finite intersections,
		\item $\bigcup_{\beta<\alpha} F_\beta\sub F_\alpha,$
		\item $|F_\alpha|<\fd$,
		\item $\sset{d_\beta, a_\beta}{ \beta<\alpha}\leq_{F_\alpha}a_\alpha.$ 
	\end{enumerate}
	Let $F_0:=\{\omega\}$ and $a_0\in A$.
	Fix $\alpha<\fd$ and assume that the elements $a_\beta\in A$ and the sets $F_\beta$, for $\beta<\alpha$, are defined. 
	For functions $s, f\in\roth$ let $s\circ f$ be the function such that $(s\circ f)(n):=s(f(n))$ for $n\in\omega$.
    For a set $B\sub \roth$, let 	\[
	\maxfin B:=\sset{\max E}{E\text{ is a finite subset of }B}\sub\roth,
	\]
where $\max E(n)=\max\{e(n):e\in E\}$ for all $n\in\w$. 
	Let 
	\[
	F:=\bigcup_{\beta<\alpha}F_\beta,\quad  
	S:= \maxfin\sset{d_\beta, a_\beta}{\beta<\alpha},\quad\text{and}\quad
	S\circ F:=\{s\circ f: s\in S, f\in F\}.
	\]
	
	\bclm
	There is $a_\alpha\in A$ such that $S\circ F\leq^{\infty}a_\alpha$.
	\label{dom}
	\eclm
	
	\bpf
	Since the set $A$ is $\fd$-unbounded, for every $y\in \roth$ the set $\sset{a\in A}{a\les y}$ has size smaller than $\fd$.
	We have $\card{S}<\fd$ and by regularity of $\fd$, we have $\card{F}<\fd$.
	Thus, $\card{S\circ F}<\fd$.
	It follows that the set 
	\[
	\bigcup_{y\in S\circ F}\sset{a\in A}{a \les y}
	\] 
	has size smaller than $\fd$.
	Therefore, there is a function
	\[
	a_\alpha\in 
	A\sm \bigcup_{y\in S\circ F}\sset{a\in A}{a \les y},
	\]
	and thus
	\[
	S\circ F\leq^\infty a_\alpha.\qedhere
	\]
	\epf
	
	\bclm
	Every intersection of finitely many elements of the set $F\cup\{[s\leq a_\alpha]: s\in S\}$ is infinite.
	\eclm
	\bpf
	The set $F$ has the finite intersection property and its elements are infinite, so it is enough to show, that for a fixed $f\in F$ and a finite subset $S'$ of the set $S$ the set
	\[
	f\cap\bigcap_{s\in S'}[s\leq a_\alpha]
	\]
	is infinite.
	Fix $f$ and $S'$ as above.
	We have
	\[
	\bigcap_{s\in S'}[s\leq a_\alpha]=\bigcap_{s\in S'}\{n:s(n)\leq a_\alpha(n)\} \supseteq [\max S' \leq a_\alpha].
	\]
	We have $\max S' \in S$, and thus $(\max S' )\circ f \leq^\infty a_\alpha$ by the definition of $a_\alpha$.
	It follows that 
	\[
	(\max S')(f(n))\leq a_\alpha(n)\leq a_\alpha(f(n))
	\]
	for infinitely many $n$.
	Thus, the set $f\cap\bigcap_{s\in S'}[s<a_\alpha]=f\cap[\max S' \leq a_\alpha]$ is infinite.
	\epf
	Closing the set $F\cup \sset{[s\leq a_\alpha]}{ s\in S}$ with respect to finite intersections of its elements, we get a set $F_\alpha$.
	Since $[s\leq a_\alpha]\in F_\alpha$ for all $s\in S$, we have 
	\[
	\sset{d_\beta, y_\beta}{ \beta<\alpha}\leq_{F_\alpha} a_\alpha.
	\]
	It finishes the construction.
	
	Let $U$ be an ultrafilter containing the family $\bigcup_{\alpha<\mathfrak{d}}F_\alpha$.
    Then the set $\sset{a_\alpha}{\alpha<\fd}$ is $\le_U$-unbounded:
	Fix $b\in\roth$.
	Since the set $D$ is dominating, there is $\alpha<\fd,$ such that $b\leq^* d_\alpha$.
	By the construction we have $b\leq^* d_\alpha\leq_U a_{\alpha+1}$, and thus $b\leq_U a_{\alpha+1}$.
    We have $\fb(\cU)=\fd$: Let $Y=\{y_\alpha:\alpha<\beta\}$ be a subset of $\roth$ for $\beta<\fd$. For every $\alpha<\beta$ we can find $d_{\alpha'}\in D$ such that $y_\alpha\leq^*d_{\alpha'}$. Since $\beta<\fd$ and $\fd$ is regular, there exists $\gamma<\fd$ such that $\{d_{\alpha'}:\alpha<\beta\}\leq_{F_\gamma}a_\gamma$ and consequently $\{y_\alpha:\alpha<\beta\}\leq_U a_\gamma$.   
	\epf
	For an ultrafilter $U$, let $2U$ be the ultrafilter generated by the set $\sset{2u}{u\in U}$, where $2u=\{2k:k\in u\}$.
	\bthm
	Assume that $X\sub\roth$ is $\leU$-unbounded for some ultrafilter $U$ with $\bof(U)=\fd$. 
	Then there is a $2U$-scale $Y$ such that  $X\times (Y\cup \Fin)$ is not Menger.
	\label{ultra}
	\ethm
	
	For functions $a, b\in\roth$ define a set
	\[
	a\uplus b := \{2k: k\in a\}\cup\{2k+1: k\in b\}\text{ for }a,b\in[\omega]^\omega.
	\]

	\blem[{\cite[Lemma 2.6]{Szewczak2017}}]
	\mbox{}
	\begin{enumerate}
		\item For each set $a\in\roth$ and each natural number $n$, we have $(a\uplus a)(2n)=2a(n)$.
		\item For all sets $a, b, c, d\in\roth$ with $a\leq b, c\leq d$, we have  $(a\uplus c)\leq (b\uplus d)$. 
	\end{enumerate}
	\label{uplus}
	\elem

	\bpf[{Proof of Theorem~\ref{ultra}}] 
    Every $\leU$-unbounded set contains a $U$-scale of the form $\sset{a_\alpha}{\alpha<\fb(U)}$, where $a_\alpha\leU\ a_\beta$ for all $\alpha\leq\beta<\fb(U)$.
	Let $A\sub X$ be such a $U$-scale.
    Since $\fb(U)=\fd$, we have $|A|=\fd$.
	Then the set 
	\[
	A_\uplus:=\sset{a_\alpha\uplus a_\alpha}{\alpha<\fb(U)}
	\]
	is a $2U$-scale:
	Fix $b\in \roth$ and let $b':=\sset{b(2n)}{n\in\w}$.
	Since $A$ is a $U$-scale, there is $\alpha<\bof(U)$ with $b'\leU a_\alpha$.
	For $n\in [b'\leq a_\alpha]$ we have
	\[
	b(2n)=b'(n)\leq a_\alpha(n)\leq 2a_\alpha(n)=(a_\alpha\uplus a_\alpha)(2n),
	\]
	and thus $2[b'\leq a_\alpha]\sub [b\leq a_\alpha\uplus a_\alpha]$.
	Since $[b'\leq a_\alpha]\in U$, we have $b\le_{2U} a_\alpha\uplus a_\alpha$.
	It follows that for each $\beta$ with $\alpha\leq\beta<\fd$, we have 
	\[
	b\leq_{2U} (a_\alpha\uplus a_\alpha)\leq_{2U} (a_\beta\uplus a_\beta).
	\]
	
	Let $D=\sset{d_\alpha}{\alpha<\fd}$ be a dominating set.
	Fix $\alpha<\bof(2U)$ and  let $d_\alpha'\in\roth$ be such that $a_\alpha,d_\alpha\leq d_\alpha'$.
	By Lemma~\ref{uplus}, we have $a_\alpha\uplus a_\alpha\leq a_\alpha\uplus d_\alpha'$.
	Then the set 
	\[
	Y:=\sset{a_\alpha\uplus d'_\alpha}{ \alpha<\fb(2U)}
	\]
	is a $2U$-scale.
	
	The space  $\PN$ with the symmetric difference operation $\oplus$ is a topological group.
	
	\bclm[{\cite[Claim 2.9]{Szewczak2017}}]
	\label{oplus}
	The set 
	\[
	(2X)\oplus (Y\cup\Fin)=\{2x\oplus y: x\in X, y\in (Y\cup\Fin)\}\sub\roth
	\]
	is dominating in $\roth$.
	\eclm
	\bpf
	For each $\alpha<\fd$ we have $a_\alpha\in X$, $a_\alpha\uplus d_\alpha'\in Y$, and
	\[
	d_\alpha\leq d_\alpha'\leq 2d_\alpha'+1=2a_\alpha\oplus(2a_\alpha\cup(2d'_\alpha +1))=2a_\alpha\oplus(a_\alpha\uplus d'_\alpha)\in
	(2X)\oplus Y.
	\]
	The inclusion $(2X)\oplus (Y\cup\Fin)\sub\roth$ follows from the definition of $\uplus$, and by the above the set $(2X)\oplus (Y\cup\Fin)$ is dominating.
	\epf
	Claim \ref{oplus} shows that there exists a continuous dominating image of the space $X\times (Y\cup\Fin)$ in $\roth$.
	Every continuous function into $\roth$ is cusco.
    By Theorem \ref{thm:Hur}, $X\times (Y\cup\Fin)$ is not Menger.
	\epf
	
	\bcor
	Assume that $\cf \fd=\fd$.
	Then for each each $\fd$-unbounded set $X$, there is a $U$-scale $Y$ for some ultrafilter $U$ such that $X\x (Y\cup\Fin)$ is not Menger. \hfill $\square$
	\label{cor:bez}
	\ecor
	
	Now we present a direct consequence of the above result.	
	A set of reals is a \emph{$\mathfrak{d}$-Lusin} set if it is homeomorphic to a subset of $\roth$ of size at least $\fd$ and whose intersection with any meager set in $\roth$ has size smaller than $\mathfrak{d}$.
	A $\mathfrak{d}$-Lusin set exists, e.g., under the assumption $\cov(\mathcal{M})=\cof(\mathcal{M})$~\cite[Lemma 29]{Scheepers_2002}.
	Let $X\sub\roth$ be a $\fd$-Lusin set. 
	Since each bounded set in $\roth$ is meager, any bounded subset of $X$ has size smaller than $\fd$.
	It follows that $X$ is $\fd$-unbounded.
	Szewczak--Tsaban proved that for any $\fd$-Lusin set there is a Menger set of reals whose product space is not Menger~\cite[Corollary 2.11]{Szewczak2017}.
	For each $U$-scale $Y$, the set $Y\cup\Fin$ is Menger in all finite powers, and thus the following Corollary is a strengthening of the mentioned result (when $\cf \fd =\fd$).
	
	\bcor\label{cor:lusin}
	Assume that $\cf\mathfrak{d}=\mathfrak{d}$. 
	For every $\mathfrak{d}$-Lusin set $X$ there is a set of reals $Y$ that is Menger in all finite powers, such that the product space $X\x Y$ is not Menger.
	\ecor

	Corollary~\ref{cor:lusin} is also related to a result of Szewczak--Wi\'{s}niewski ~\cite[Theorem 3.1]{szewczak2019productsluzintypesetscombinatorial} that assuming $\cov(\cM)=\cof(\cM)$, for each $\fd$-Lusin set $X$ there is a $\fd$-Lusin set $Y$ that is Menger in all finite powers, such that the product space $X\x Y$ is not Menger.

	\section{Productively Scheepers and productively Hurewicz spaces}\label{sec:Result} 
	
	In this section we prove the following Theorem.
	
	\bthm\label{result}
	Assume that $\fb=\fd$.
	\be
	\item
	Let $X$ be a space such that for every ultrafilter $U$ and $U$-Menger space $T$ the product space $X\x T$ is Menger.
	Then the space $X$ is productively Hurewicz. 
	
	\item
	Each productively Scheepers space is productively Hurewicz.
	\item 
	Each space which is productively $U$-Menger for all ultrafilters $U$ is productively Hurewicz.
	\ee
	\ethm
	
	Before we prove this theorem, we need auxiliary results.
	
	\blem\label{lem:U}
	Assume that $\fb=\fd$.
	Then each unbounded set is $\leU$-unbounded for some ultrafilter $U$.
	\label{unb}
	\elem

	\bpf
	Let $X\sub\roth$ be an unbounded set.
	In light of Theorem~\ref{ultra} it is enough to prove that $X$ contains a $\fd$-unbounded set.
	Let $D=\sset{d_\alpha}{\alpha<\fd}$ be a dominating set. 
	Let $x_0\in X$.
	Fix $\alpha<\fd$ and assume that elements $x_\beta\in X$, $\beta<\alpha$, have been already defined.
	Since $\fb=\fd$, there is a function $b$ such that $\{x_\beta, d_\beta:\beta<\alpha\}\les b$.
	The set $X$ is unbounded, so there is a function $x_\alpha\in X$ such that
	\[
	\{x_\beta, d_\beta:\beta<\alpha\}\leq^*a_\alpha\leq^\infty x_\alpha.
	\]
	Then $\sset{x_\alpha}{\alpha<\fd}$ is $\fd$-unbounded.
	\epf

    \bcor
	Assume $\fb=\fd$. Then Hurewicz is equivalent to being $U$-Menger for all ultrafilters $U$.   
    \ecor
	
	\blem [{\cite[Theorem 5.2]{Szewczak2017}}]
	Let $U$ be an ultrafilter and $X\sub\roth$ be a $U$-scale.
	Then the space $X\cup\Fin$ is productively $U$-Menger.
	\label{bez}
	\elem
	
	\bpf[{Proof of Theorem~\ref{result}}]
	(1)
	Let $H$ be a Hurewicz space and suppose that $X\times H$ is not Hurewicz.
	By the assumption, the product space $X\times H$ is Menger, so in particular it is Lindel\"of.
	By Theorem~\ref{thm:Hur}, there is a cusco image $Z\sub\roth$  of  $X\times H$ which is unbounded.
	By Lemma~\ref{unb}, the set $Z$ is $\leU$-unbounded for some ultrafilter $U$.
	By Theorem~\ref{ultra}, there is a $2U$-scale $Y$, such that $Z\times (Y\cup \Fin)$ is not Menger.
	The space $Z\times (Y\cup\Fin)$ is a cusco image of the space $X\times H\times (Y\cup \Fin)$.
	Since the Menger property is preserved by cusco maps, the product $X\times H\times (Y\cup \Fin)$ is not Menger.
	
	On the other hand, since the space $H$ is Hurewicz, it is also $2U$-Menger.
	By Theorem~\ref{bez}, the product space $H\times (Y\cup \Fin)$ is $2U$-Menger.
    Since $H$ is hereditarily Lindel\"of and $Y\cup\Fin$ is a set of reals, this product is also hereditarily Lindel\"of.
	By the assumption on $X$, the product space $X\times H\times (Y\cup \Fin)$ is Menger, a contradiction.
	
	Items (2) and (3) follow directly from (1).
	\epf

	\bthm\label{thm:pMpUM}
	Each productively Menger space is productively $U$-Menger for all ultrafilters $U$ with $\bof(U)=\fd$.
	\ethm
	
	In order to prove Theorem~\ref{thm:pMpUM} we need the following observation.
	
	\blem\label{Men}
	Let $U$ be an ultrafilter. Then any $U$-Menger space is $2U$-Menger.
	\elem
	
	\bpf

    Let $X$ be a $U$-Menger set and $Y\sub\roth$ be its cusco image.
	Let $y'\in\roth$ be a function such that $y'(n):=y(2n)$ for all $n$ for each $y\in \roth$ and define $Y':=\{y':y\in Y\}$.
	Since $Y'$ is a cusco image of $X$, it is $\leU$-bounded by some function $b$.
    Fix $y'\in Y'$.
    Since $\sset{n}{y'(n)\leq b(n)}\in U$, we have
	\[
	\sset{n}{y'(n)\leq b(n)}=\sset{n}{y(2n)\leq b(n)}\sub\sset{n}{y(2n)\leq b(2n)}\in U,
	\]
	and thus $\sset{2n}{ y(2n)\leq b(2n)}\in2U$.
	We conclude that $Y$ is $\le_{2U}$-bounded by the function $b'$, hence $X$ is $2U$-Menger.
	\epf
	
	\bpf[{Proof of Theorem~\ref{thm:pMpUM}}]
	Let $X$ be a productively Menger space and suppose that there is a $U$-Menger space $H$ for some ultrafilter $U$ such that $X\times H$ is not $U$-Menger.
	By the properties of $X$ and $H$, the product $X\times H$ is Menger, in particular it is Lindel\"of.
	By Theorem~\ref{thm:Hur}, there is a cusco image $Z\sub\roth$  of  $X\times H$ which is $\leU$-unbounded.
	By Theorem \ref{ultra}, there is a $2U$-scale $Y$, such that $Z\times (Y\cup \Fin)$ is not Menger.
	The space $Z\times (Y\cup\Fin)$ is a cusco image of the space $X\times H\times (Y\cup \Fin)$.
	Since the Menger property is preserved by cusco maps, the product $X\times H\times (Y\cup \Fin)$ is not Menger.
	
	On the other hand, by Lemma~\ref{Men} and Lemma \ref{bez}, the product space $H\times (Y\cup \Fin)$ is Menger. Since $H$ is hereditarily Lindel\"of and $Y\cup\Fin$ is a set of reals, this product is also hereditarily Lindel\"of.
	By the assumption on $X$, the product space $X\times H\times (Y\cup \Fin)$ is Menger, a contradiction.
	\epf

\bthm
It is consistent with \CH{} that in the class of sets of reals there is a productively Hurewicz set which is not productively $U$-Menger for some ultrafilter $U$.
\ethm

\bpf
This example is essentially due to Repov\v{s}--Zdomskyy~\cite[Section~3]{Repov__2025}.
Let $\bbC_{\w_1}$ be an iterated Cohen forcing of length $\w_1$ with finite supports and $G$ be a $\bbC_{\w_1}$-generic over a model $V$ satisfying \CH{}.
In $V[G]$, there are a productively Hurewicz set and a Menger set from the ground model, whose product is not Menger ~\cite[Section~3]{Repov__2025}. 

We show that every ground model set $X\sub\Cantor\cap V$ is $U$-Menger in the extension, for some ultrafilter $U$.
Let $c_\alpha$ be the $\alpha$-th Cohen real added to the ground model, where $\alpha<\w_1$.
The sets 
\[
F_\alpha:=\sset{[z\leq c_\alpha]}{z\in \Un_{\beta<\alpha}\NN\cap V[\sset{c_\beta}{\beta<\alpha}]}
\]
are centered and $\Un_{\beta<\alpha}F_\beta$ is centered too, where $\alpha<\w_1$.
Let $U$ be an ultrafilter in $V[G]$ containing $\Un_{\alpha<\w_1}F_\alpha$.
Let $\varphi\colon X\to \roth$ be a continuous function in $V[G]$.
Then there is an ordinal number $\alpha<\w_1$ such that $f$ is coded in $V[\sset{c_\beta}{\beta<\alpha}]$.
Then $\varphi[X]\le_U c_\alpha$, and thus the set $X$ is $U$-Menger.
\epf
	
\section{Near coherence of filters and the Scheepers property}\label{sec:Sch-UMen}

\emph{Near coherence relation} on the set of ultrafilters 
was introduced by Blass~\cite{ncf1}.
Ultrafilters $U,U'$ are \emph{nearly coherent} if there exists a finite-to-one map $\varphi\colon\w\to\w$ such that 
$\varphi[U]=\varphi[U']$.
Following \cite{ncf1} we denote by NCF the statement that every two ultrafilters on $\w$ are nearly coherent.
NCF is consistent with ZFC~\cite{242}, it follows from $\fu<\fg$~\cite{BlaLaf89},
but is not equivalent to $\fu<\fg$~\cite{MilShe09}.
NCF yields  $\fb(U)=\fd$ for all ultrafilters $U$ ~\cite[Theorem~16]{ncf1} and it implies the existence of an ultrafilter generated by
less than $\fd$ sets.
It is known (see, e.g., the discussion in the first paragraph on \cite[P. 729]{BlaMil99}) that near coherence of two ultrafilters is witnessed by a monotone surjection $\varphi\colon\w\to\w$.
We use these results in the following Proposition.

\begin{prop} \label{expl_nachher}
Assume that NCF holds.
Every Scheepers space is $U$-Menger for all ultrafilters $U$.
Consequently, the properties $U$-Menger are equivalent for all ultrafilters $U$.
\end{prop}

\begin{proof}
Apply Theorem~\ref{thm:Hur}.
Let $Y\sub \roth$ be a cusco image of a Scheepers space $X$.
There is a function 
$h\in \roth$ witnessing that 
$Y$ is not finitely-dominating, which implies that there exists an
ultrafilter $U'$ containing $\sset{[y\leq h]}{y\in Y}$.
By NCF, there is an ultrafilter $W$ with a base of size $\kappa<\fd$ and there is a monotone surjection $\varphi\colon\w\to\w$   
with $\varphi[U']=\varphi[W]$.
Let $\sset{w_\alpha}{\alpha<\kappa}$ be a base for $W$.
The sets
\[
K_\alpha:=\sset{z\in \roth}{\forall{n\in \varphi[w_\alpha]}\exists{k\in \varphi^{-1}(n)}\, (z(k)\leq h(k)) }
\]
are compact for all $\alpha<\kappa$.
For every $y\in Y$ there exists $\alpha<\kappa$ such that 
$\varphi[w_\alpha]\sub\varphi\bigl[[y\leq h]\bigr]$, and thus $y\in K_\alpha$.
It means that $Y\sub\Un_{\alpha<\kappa}K_\alpha$.

Let $U$ be an arbitrary  ultrafilter.
Since $\fb(U)=\fd$ and each compact subspace of $[\w]^\w$ is dominated by a single
element, the set $\bigcup_{\alpha\in\kappa}K_\alpha$ is $\leq_U$-bounded,
and thus the set $Y$ is $\leq_U$-bounded, too.
\end{proof}

\section{Comments and open problems}\label{sec:q}

Under the stronger assumption $\fd=\aleph_1$, the main results established above for hereditarily Lindel\"of spaces extend to arbitrary topological spaces.

\bthm
Assume $\fd=\aleph_1$.
In the class of general spaces the following assertions hold.
\begin{enumerate}
		
		\item Each productively Menger space is productively $U$-Menger for all ultrafilters $U$.
		
		\item Each space which is productively $U$-Menger for all ultrafilters $U$ is productively Hurewicz.
		
		\item Each productively Scheepers space is productively Hurewicz.
		
	\end{enumerate}
	\ethm
	\bpf
    Since $\fd=\aleph_1$, for every ultrafilter $U$, we have $\bof(U)=\aleph_1$.
    Then, by the result of Szewczak--Tsaban~\cite[Theorem~3.10(4)]{szewczak2017productsgeneralmengerspaces}, Lemma~\ref{bez} holds in the class of general spaces. 
	In~(1), we proceed in a similar manner as in the proof of Theorem~\ref{thm:pMpUM}.
	In~(2) and~(3), we proceed in a similar manner as in the proof of Theorem~\ref{result}.
	\epf
	
	The following questions arise naturally from our investigations.
	
	\bprb
	Is it consistent that:
	\be 
	\item There is a space that is productively $U$-Menger for all ultrafilters $U$, but not productively Scheepers? 
    \item There is a productively Scheepers space that is not productively $U$-Menger for some ultrafilter $U$? 
	\item There is a productively Menger space that is not productively Scheepers?  
	\item There is a productively Scheepers space that is not productively Menger?
	\ee 
	\eprb
	
	Repov\v{s}--Zdomskyy~\cite[Theorem~1.1]{Repov__2025} proved that, in the Laver model (where $\fb=\fd$), any Hurewicz set of reals is productively Menger and also productively Hurewicz.
	However, their result does not capture the Scheepers property.
	\bprb
	In the Laver model, is every Hurewicz space productively Scheepers?
	\eprb

\section*{Acknowledgements}
The result that, under $\fb=\fd$, every productively Scheepers space is productively Hurewicz was established in the Master’s thesis of the first-named author~\cite{master}.
We are grateful to Professor Piotr Zakrzewski, who reviewed the thesis, for carefully reading the manuscript and sharing his comments with us.

We are also grateful to the reviewer for the corrections and suggestions that substantially improved the quality of the paper.

\bibliographystyle{abbrv}
	\bibliography{biblio}
	
\end{document}